\documentclass[11pt]{article}

\input{amssym.def}
\input{amssym}
\usepackage{eufrak}

\newtheorem{theorem}{Theorem}

\newtheorem{proposition}[theorem]{Proposition}
\newtheorem{remark}[theorem]{Remark}

\def\hl{\hat{l}}
\def\hL{\hat{L}}
\def\hhh{{\frak{h}}}

\def\hd{\hat{d}}

\def\qq{q^{-1}}
\def\H{\Bbb{H}}

\def\Tr{\mathrm{Tr}}
\def\Trr{\Tr_R}
\def\SS{\Bbb{S}}
\def\PP{\cal{P}}

\def\hLL{\hat{{\cal L}}}
\def\LL{{\cal{L}}}

\def\PP{{\cal P}}
\def\gg{{\frak g}}

\def\si{{\sigma}}

\def\de{\delta}

\def\ot{\otimes}
\def\C{{\Bbb C}}

\def\vv{V^{\otimes 2}}
\def\RR{R^{-1}}
\def\AA{{\cal A}}

\def\ov{\overline}
\def\hmu{\hat{\mu}}
\def\hnu{\hat{\hnu}}
\def\hd{\hat{d}}
\def\hf{\hat{f}}
\def\la{{\lambda}}
\def\vl{{V_\lambda}}

\def\be{\begin{equation}}
\def\ee{\end{equation}}

\begin{document}

\makeatletter
\renewcommand{\theequation}{{\thesection}.{\arabic{equation}}}
\@addtoreset{equation}{section} \makeatother

\title{Generalized Harish-Chandra morphism on Reflection Equation algebras}
\author{\rule{0pt}{7mm} Dimitry
Gurevich\thanks{gurevich@ihes.fr}\\
{\small\it Institute for Information Transmission Problems}\\
{\small\it Bolshoy Karetny per. 19,  Moscow 127051, Russian Federation}\\
\rule{0pt}{7mm} Pavel Saponov\thanks{Pavel.Saponov@ihep.ru}\\
{\small\it
National Research University Higher School of Economics,}\\
{\small\it 20 Myasnitskaya Ulitsa, Moscow 101000, Russian Federation}\\
{\small \it and}\\
{\small \it
Institute for High Energy Physics, NRC "Kurchatov Institute"}\\
{\small \it Protvino 142281, Russian Federation}}

\maketitle

\begin{abstract}
We consider the so-called generalized Harish-Chandra morphism, taking the center of the enveloping algebra $U(gl(N))$ to the commutative algebra generated by eigenvalues of the
generating matrix of this algebra, and generalize this construction to Reflection Equation algebras. To this end we introduce the eigenvalues of the generating matrix of the Reflection Equation algebra (modified or not), corresponding  to a skew-invertible Hecke symmetry and define  the  generalized Harish-Chandra morphism in a similar way. We use this map in 
order to introduce quantum analogs of the so-called weight systems.
\end{abstract}

{\bf AMS Mathematics Subject Classification, 2010:} 17B37,  81R50


\section{Introduction}
Consider a Lie algebra $gl(N)$ and fix in it the standard basis $\hl_i^j$ such that
$$
\hl_i^j\, \hl_k^m-l_k^m\, \hl_i^j=\hl_i^m\, \de_k^j- \hl_k^j\, \de_i^m,\quad  1\leq i,j,k,m \leq N.
$$
The matrix $\hL=\|\hl_i^j\|_{1\leq i, j \leq N}$\footnote{The lower index of entries numerates the rows of the corresponding matrix.}
meets a polynomial relation $\PP(\hL)=0$, which is called the Cayley-Hamilton identity. The coefficients of the corresponding polynomial
\be
\PP(t)=\sum_{k=0}^N (-1)^k a_k\, t^{N-k}
\label{CH}
\ee
belong to the  center $Z(U(gl(N))$ of the universal enveloping algebra $U(gl(N))$. Besides, the polynomial $\PP(t)$ is unital, i.e. $a_0=1$.

It is natural to introduce the roots $\hmu_i$, $1\leq i\leq N$ of the polynomial $\PP$ as follows:
\be
\sum_{i=1}^N\, \hmu_i=a_1 \qquad \sum_{i<j}^N\, \hmu_i\, \hmu_j=a_2\quad \dots\quad \prod_{i=1}^N\, \hmu_i=a_N.
\label{sta}
\ee
These roots are assumed to be central in the extended algebra $U(gl(N))[\hmu_1,\dots,\hmu_N]$ and are called  the {\em eigenvalues} of the matrix $\hL$.
Moreover, any central element of the algebra $U(gl(N))$ can be expressed as a symmetric polynomial in the eigenvalues $\hmu_i$.
Thus, we get an isomorphism
\be
\varphi: Z(U(gl(N)))\to \mathrm{Sym}(\hmu), \label{HC}
\ee
where $\mathrm{Sym}(\hmu)$ stands for the space of all symmetric polynomials in $ \hmu_i$. We call this map the  {\em generalized  Harish-Chandra morphism}.

Note that usually this map is called the Harish-Chandra  isomorphism, it maps the center $Z(U(gl(N)))$ to the commutative algebra of "shifted symmetric polynomials". However,
these polynomials become symmetric after a linear change of the generators.

It is known that any finite-dimensional irreducible $U(gl(N))$-modules $\vl$ is labeled by a partition
\be
\la=(\la_1\geq \la_2\geq,\dots ,\geq \la_N), \label{part}
\ee
where $\la_i$ are non-negative integers.

Besides, any central element $\cal{ B}$ being represented in such a $U(gl(N))$-module becomes a scalar operator ${\chi_\la(\cal B)}\, I$, where $\chi_\la(\cal B)$   is a  number,
called the {\em character} of $\cal B$. The characters of the elements ${\cal B}=Tr \hL^k$ have been computed in \cite{PP}. The result of \cite{PP} agrees with the characters
of the eigenvalues $\hmu_i$, defined by the formula
$$
\chi_\la(\hmu_i)=\la_i+N-i
$$
for an appropriate ordering of $\hmu_i$.

The  purpose of the current paper is to define the  quantum analog of the Harish-Chandra morphism on  the Reflection Equation (RE) algebras, associated to  Hecke symmetries
of general type (see the next section) and to exhibit an  application of this quantum Harish-Chandra morphism.
The point is that the generating  matrix of any such  algebra satisfies a Cayley-Hamilton relation
with central coefficients\footnote{Observe that this property for a matrix with noncommutative entries is somewhat
exceptional. Besides the enveloping algebras of some Lie algebras (see, for instance \cite{Go}) and the RE algebras, we can mention only  the quantum
groups $U_q(gl(N))$, which enable us to introduce such  matrices.  However, it is so in virtue of  their relations with the RE algebras, associated with the
Drinfeld-Jimbo Hecke symmetries (see the last section).}.

This property enables us to introduce the eigenvalues of the generating matrices and to parameterize (hopefully) all central elements of the corresponding algebra by these
eigenvalues. The parametrization of some central elements (the power sums and Schur polynomials) has been found  in \cite{GPS1, GPS2, GPS4}. Namely, by associating
such a parameterized form with central elements we get a map which is called {\em the quantum Harish-Chandra morphism}.

As for the aforementioned application, we introduce quantum analogs of the so-called weight systems. The modern viewpoint onto the usual weight system is exhibited in \cite{Y1, Y2}.
With any element of the symmetric group the author associates a central element of an algebra $U(gl(m))$ or $U(gl(m|n)$. On the second step by applying the Harish-Chandra morphism,
the center of the enveloping algebra is sent into the enveloping algebra of the Cartan subalgebra.

Our quantum analog of this scheme  is also composed of two maps. The first map sends the  Hecke algebra to the center of the RE algebra. The second map is the
quantum Harish-Chandra morphism\footnote{However,  so far it is not clear whether such quantum analog meets the four terms relation (see \cite{Y1}). We plan to consider this
problem in the future.}. Note that for the sake of simplicity in the bulk of the paper we consider the Hecke symmetries of $GL(N)$ type also called  {\it even symmetries}.

The case of the Hecke symmetries of general type is considered in the last section. Also, there we consider the algebras related to the Drinfeld-Jimbo
symmetries. Note that in this case the construction of the irreducible representations of the RE algebras and computation of  the characters of the center can be done by the methods
of  \cite{RTF}. Note that exactly these methods were used for computing the values of the ``quantum Gelfand invariants'' in \cite{JLM}.

The paper is organized as follows. In the next section we recall some elements of the so-called $R$-matrix technique. In Section 3 we reproduce formulae expressing the power sums
via the eigenvalues of the generating matrices of the RE algebras and their modified counterparts under assumption that the initial Hecke symmetry is even.  In Section 4 we introduce two quantum versions (related to modified and
non-modified RE algebras) of the weight systems and exhibit a few examples. In Section 5  we briefly consider the case of general type Hecke symmetries and
compare our method and that  of \cite{JLM}, where similar results were obtained in a particular case related to the quantum group $U_q(gl(N))$.

{\bf Acknowledgement} We are very indebted to Alexander Molev for fruitful discussions.

 \section{Preliminaries}

 The ground field is assumed to be ${\Bbb C}$ and $q$ is assumed to be generic, i.e. such that $k_q\not=0$ for any integer $k$ . Hereafter, we use  the standard notation
$k_q = \frac{q^k - q^{-k}}{q-q^{-1}}$.

 Let $V = {\Bbb C}^N$ and $R\in \mathrm{End}(V^{\otimes 2})$. An operator $R$ is called {\it a Hecke symmetry} if it satisfies the the braid relation
$$
R_{12}\, R_{23}\,R_{12}= R_{23}\,R_{12}\, R_{23}
$$
(such operators are called braidings), and, in addition, meets the Hecke condition
$$
(q\, I-R)(\qq\, I+R)=0,\quad q\not= \pm 1.
$$
In what follows the symbol $I$ stands for an identity operator or matrix of an appropriate size. The braid relation is written for operators $R_{12} = R\otimes I$ and
$R_{23} = I\otimes R$ from $\mathrm{End}(V^{\otimes 3})$. Below we use the following embeddings of operator $R$ in the larger spaces
$$
R_{kk+1} = I^{\otimes (k-1)}\otimes R\otimes I^{\otimes{p-k-1}}\qquad 1\le k\le p-1.
$$
Here $p\geq 2$ is an integer. In what follows we do not fix a specific value of $p$, just assuming it to be large enough for all matrix formulas to make sense. Also, we use the 
shorthand notation $R_k = R_{kk+1}$. In a similar way, for $\forall\, X\in \mathrm{Mat}_N$ we set
$$
X_k = I^{\otimes (k-1)}\otimes X\otimes I^{\otimes{p-k}}\quad 1\le k\le p.
$$

There are known many examples of Hecke symmetries which are deformations of the usual flip $P$ or super-flips $P_{m|n}$ and big families of Hecke symmetries, which are
deformations neither of the flips nor of the super-flips. In the classical case (i.e. if $R=P$) the Hilbert-Poincar\'e series of the skew-symmetric algebra $\Lambda(V)$ of the space
$V$ is the polynomial $(1+t)^N$ of degree $N$. In general it is not so even for the Hecke symmetries of $GL(N)$ type. By such a symmetry  we mean $R$, for which the Hilbert-Poincar\'e series of the corresponding $R$-skew-symmetric algebra
\be
\Lambda_R(V)=T(V)/\langle \mathrm{Im}(\qq\, I+R)\rangle,\qquad T(V)=\bigoplus_{k\ge 0} V^{\ot k}
\label{vn-alg}
\ee
of the basic space $V$ is a  polynomial. Here, $\langle J \rangle$ is the ideal generated in the free tensor algebra $T(V)$ by a subset $J$.
If the degree of this polynomial is $m$, we say that the bi-rank of $R$ is $(m|0)$. In general, $m\leq N$. Also, we say that the bi-rank of a given Hecke symmetry $R$ is $(m|n)$ if
the Hilbert-Poincar\'e series of the corresponding $R$-skew-symmetric algebra $\Lambda (V)$ is a rational function with the numerator of degree $m$ and the denominator of degree
$n$.

It is worth noticing that any Hecke symmetry defines a so-called $R$-matrix representation of a Hecke algebra. Let us remind that the $A_{n-1}$ type Hecke algebra $\H_n(q)$ is a
unital associative algebra over the complex field $\Bbb{C}$ generated by the set of  generators $\tau_i$, $1\le i\le n-1$,  satisfying the following relations:
\begin{eqnarray}
\tau_i\, \tau_{i+1}\, \tau_i=\tau_{i+1}\, \tau_i\, \tau_{i+1},&\quad& 1\le i \le n-2,\nonumber\\
\rule{0pt}{4mm}
\tau_i\, \tau_j=\tau_j\, \tau_i\, &\quad& |i-j|\geq 2,\label{H-alg}\\
\rule{0pt}{4mm}
\tau_i^2=e+(q-q^{-1})\tau_i,&\quad& 1\le  i\le n-1,\nonumber
\end{eqnarray}
where $e$ is a unit element of the algebra and $q\in\Bbb{C}\setminus \{0,\pm 1\}$ is a numeric parameter. The limit of these algebras (which are naturally embedded in each
other) is denoted by $\H(q)$ and is also called the Hecke algebra. In the same manner we denote by $\SS$ the limit of the symmetric groups $\SS_n$.

Note that at $q=\pm 1$ the defining relations on the generators $\tau_i$ turn into those for the generators of the symmetric group $\SS_n$. For a generic value of $q$ the Hecke
algebra $\H_n(q)$ is a deformation of the group algebra $\C[\SS_n]$.

The aforementioned $R$-matrix representation $\rho_R$ of the Hecke algebra $\H_n(q)$ is defined on generators $\tau_i$ as follows:
$$
\rho_R: \H_n(q)\rightarrow \mathrm{End}(V^{\otimes n})\qquad
\rho_R(\tau_i)=R_i,\quad 1\le i\le n-1.
$$

We assume all Hecke symmetries, we are dealing with, to be skew-invertible. This means that there exists an operator $\Psi:\vv\to \vv$ such that
\be
\Tr_{(2)} R_{12}\, \Psi_{23}=P_{13}\quad \Longleftrightarrow \quad \sum_{ia,b =1}^NR_{ib}^{ja} \Psi_{ak}^{b\,n}=\de_i^n\de_k^j.
\label{Psi}
\ee
Hereafter, we assume a basis of the space $V$ to be fixed and all operators to be represented by their matrices in this basis.

Also, we need the $N\times N $  matrix $C=\|C_i^j\|$, where
$$
C_i^j=\sum_{k=1}^N \Psi_{ik}^{jk}.
$$
This matrix is used in the definition of the $R$-trace: for an $N\times N$ matrix $A$ we define its $R$-trace by the rule
$$
\Tr_R A= \Tr (C A).
$$

Now, introduce two forms of the RE algebras. Given a skew-invertible Hecke symmetry $R$,  the {\it RE algebra} $\LL(R)$ is a unital associative algebra generated by entries of the
matrix $L=\|l_i^j\|_{1\leq i,j \leq N}$ subject to homogenous matrix relation:
\be
R\, (L\ot I)\, R\, (L\ot I)-(L\ot I)\, R\, (L\ot I)\, R = 0.
\label{REA}
\ee

The {\it modified RE algebra} $\hLL(R)$ is generated  by entries of the matrix $\hL=\|\hl_i^j\|_{1 \leq i, j \leq N}$ subject to the relation:
\be
R\, (\hL\ot I)\, R\, (\hL\ot I)-(\hL\ot I)\, R\, (\hL\ot I)\, R=R\, (\hL\ot I)-(\hL\ot I)\, R.
\label{odin}
\ee

Note that the algebras $\hLL(R)$ and $\LL(R)$ are isomorphic to each other provided $q\not = \pm 1$. Their isomorphism is generated by the following relation between the matrices $\hL$ and $L$
\be
L=I - (q-\qq)\, \hL\quad \Longleftrightarrow \quad l_i^j= \delta_i^j - (q-\qq)\hl_i^j.
\label{sh}
\ee

Now, introduce the following notations
$$
L_{\ov 1}=L_{ 1},\quad L_{\ov k}=R_{k-1}\,L_{\ov{k-1}}\,\RR_{k-1}, \quad  k\ge 2
$$
and
$$
L_{\ov{1\to k}}=L_{\ov 1}\, L_{\ov 2}\dots L_{\ov k}.
$$

As was shown in \cite{IP}, for any element $z\in \H_n(q)$ the following homogeneous polynomial in the generators of the RE algebra
$$
ch_n(z)=\Tr_{R(1\dots n)}\left(\rho_R(z) L_{\ov{1\to n}}\,\right) = \Tr_{R(1\dots n)} \left( L_{\ov{1\to n}}\,\,\rho_R(z)\right)
$$
is central in $\LL(R)$.

Thus, we have a set of maps
$$
ch_n: \H_n(q)\to Z(\LL(R)),\quad z\mapsto ch_n(z),\quad \forall\, n\ge 1,
$$
and consequently  the map
\be
ch: \H(q)\to Z(\LL(R)). \label{elem}
\ee

The linear envelope of the elements (\ref{elem}) for all $z\in \H(q)$  forms a {\em characteristic subalgebra}\footnote{Conjecturally, for any Hecke symmetry $R$
the characteristic subalgebra coincides with the whole center $Z(\LL(R))$. \label{fn}} ${\cal C}(\LL(R))\subset Z(\LL(R))$ of the center  $Z(\LL(R))$.

Some particular elements $z\in \H(q)$ map into elements $ch(z)\in Z(\LL(R))$ of special interest. Namely, they are the so-called quantum Schur polynomials and power sums.
The Schur polynomials $s_\la(L)$ are labeled by partitions (\ref{part}). For an explicit description of the Schur polynomials the reader is referred to \cite{GS2}.
We only note that they generate  an associative algebra since
\be
s_\la(L) \, s_\mu(L)= {\sum}_\nu C_{\la\, \mu}^\nu s_\nu(L)
\label{Lit-Rich}
\ee
with the classical Littlewood-Richardson coefficients $C_{\la\, \mu}^\nu $. The full set of quantum Schur polynomials $s_\lambda(L)$  forms a linear
basis of the characteristic subalgebra  ${\cal C}(\LL(R))$.

Below, we need a particular subclass of quantum elementary symmetric polynomials. They corresponds to partitions $\lambda =(1^k)$. First, we  define the set of
skew-symmetrizers $\AA^{(k)}$ by the following rules:
$$
\AA^{(1)}=I,\quad \AA^{(k)}=\frac{1}{k_q} \AA^{(k-1)}\left(q^{k-1} I-(k-1)_q \, R_{k-1}\right)\AA^{(k-1)}.
$$
Note that $\AA^{(k)}$ maps $V^{\otimes k}$ into $\Lambda_R^{(k)}(V)$ which is the $k$-th homogenous component of the $R$-skew-symmetric algebra  (\ref{vn-alg}) of
the space $V$.

The elementary symmetric polynomials $e_k(L)$ are then defined by the $R$-trace:
$$
e_k(L) = \Tr_{R(1\dots k)}(\AA^{(k)}L_{\overline{1\rightarrow k}}).
$$

If in (\ref{elem}) we put $z_{k}=\tau_{k-1}\, \tau_{k-1}\dots \tau_{1}$, we get a central element $p_k(L)$, called the $k$-th power sum:
$$
p_k(L) = \Tr_{R(1\dots k)}(R_{k-1}\dots R_1\,L_{\overline{1\rightarrow k}})
$$
Note that it can be simplified to $p_k(L)=\Trr\, L^k$. In more detail we consider the power sums in the next section.

The generating matrix $L$ of the RE algebra $\LL(R)$ satisfies the Cayley-Hamilton identity:
\be
L^m-q\, e_1(L)\, L^{m-1}+q^2\, e_2(L)\, L^{m-2}+...+ (-q)^{m-1} \, e_{m-1}(L)L+(-q)^{m}  e_{m}(L)\, I=0,
\label{CH1}
\ee
where $(m|0)$ is the bi-rank of $R$.

Now we introduce the quantum eigenvalues $\{\mu_i\}_{1\le i\le m}$ as the roots of this polynomial:
\be
\sum_{i=1}^m \, \mu_i=q\, e_1(L) , \qquad \sum_{i<j}^m\, \mu_i\, \mu_j=q^2e_2(L)\quad \dots \quad \prod_{i=1}^m \mu_i=q^{m}\, e_m(L).
\label{q-mu}
\ee

Since all elements $e_k(L)$ are central, it is natural to assume the eigenvalues $\mu_i$ to be central in the extended algebra
$\LL(R)[\mu_1,\dots ,\mu_m]$. With the use of (\ref{q-mu}) we can rewrite the identity (\ref{CH1}) in a factorized form:
$$
\prod_{i=1}^m (L-\mu_i\, I)=0.
$$

Due to (\ref{Lit-Rich}) the elements $e_k(L)$, $1\le k\le m$, generate the whole characteristic subalgebra ${\cal C}(\LL(R))$. So, the parameterization (\ref{q-mu})
defines a map of ${\cal C}(\LL(R))$ into the ring of symmetric polynomials in the eigenvalues $\mu_i$:
\be
\varphi_q: {\cal C}(\LL(R)) \rightarrow \mathrm{Sym}(\mu).
\ee
Up to the conjecture mentioned in the footnote \ref{fn}, we get a quantum analog of the map (\ref{HC}). Namely, the map $\varphi_q$ is called
{\it the quantum  Harish-Chandra morphism}.

\begin{remark}\rm Note that in \cite{GS1} we have found the Cayley-Hamilton polynomial ${\cal Q}(t)$ for the matrix $\hL$. It reads:
$$
{\cal Q}(t) = \Tr_{R(1\dots m)}\left(\AA^{(m)}\, \prod_{k=1}^m \,\left(q^{2(k-1)}(t-q^{-k+1}(k-1)_q)I-\hL_{\ov k}\right) \right).
$$
\end{remark}

\section{Power sums and their representation in modules $\vl$}

First, we briefly remind a method of constructing the modules $\vl$.

Initially, the representation category of the modified RE algebras $\hLL(R)$ was constructed in \cite{GPS3}. Let $\{x_i\}_{1\le i\le N}$ be a fixed basis in the space $V$.
Then, as it was shown in  \cite{GPS3}, the linear action of $\hLL(R)$ generators
\be
\hl_i^j \triangleright x_k= B_k^j\, x_i
\label{act}
\ee
defines a representation of the modified RE algebra $\hLL(R)$ in the space $V$. The $N\times N$ matrix $B = \|B_i^j\|$ is of the form:
$$
B_i^j=\sum_{k=1}^N\Psi_{ki}^{kj}
$$
where $\Psi$ is the skew-inverse to $R$ (\ref{Psi}).

Besides,  the action (\ref{act}) of the algebra $\hLL(R)$ was extended on the tensor powers $V^{\ot k}$, $k\ge 2$. These  $\hLL(R)$-modules were decomposed in the direct sums 
of the invariant submodules\footnote{The irreducibility of $V_\lambda$ is an open problem.} $\vl$. The modules $V_\lambda$ are images of projection operators $P_\lambda(R)$ 
which are used in the construction of the quantum Schur polynomials $s_\lambda(L)$.

Lately, the same results for the representation theory of the RE algebra $\LL(R)$ was obtained in terms of the quantum double construction. An example of such a quantum double is
the pair  $(\LL(R), \,\bigoplus_k V^{\ot k})$. The $\LL(R)$-modules $\vl$ can be cut out from $T(V)$ by means of the same projectors $P_\lambda$.

\begin{remark}\rm
We can also construct the representations of the RE algebra in the dual spaces $(V^*)^{\otimes k}$ and extract invariant submodules $V^*_\lambda$ in a similar way.
We do not need them below.
\end{remark}

Though the irreducibility of  the modules $\vl$ is not proved for general Hecke symmetry $R$, the images of all central elements of $\LL(R)$ or $\hLL(R)$ turn out to be scalar operators.
In order to compute the character of any central element $z\in Z(\LL(R))$ (or $z\in Z(\hLL(R))$) on the module $V_\lambda$ it suffices to express $z$ in terms of the quantum eigenvalues
$\mu_i$ (respectively,  $\hmu_i$) and use the following result from \cite{GSZ}.

\begin{proposition}
Let $R$ be skew-invertible Hecke symmetry of bi-rank $(m|0)$ and $\mu_i$, $1\leq i\leq m$ be the quantum eigenvalues of the generating matrix $L$ of the corresponding RE algebra 
$\LL(R)$. Then in any module $V_\lambda$ the eigenvalues $\mu_i$ are represented by scalar operators. Moreover, there exists an ordering of the eigenvalues $\mu_i$ such that the
corresponding scalar operators read:\rm
\be
\mu_i = \mu_i(\lambda) \mathrm{Id}_{V_\lambda},\qquad
\mu_i(\la)= q^{-2(\la_i+m-i)}.
\label{mul}
\ee
\end{proposition}

In virtue of (\ref{sh}) the eigenvalues of the matrices $L$ and $\hL$ are related as follows:
\be
\mu_i=1-(q-q^{-1})\,\hmu_i. \label{zamena} \ee
This relation enables us to compute the quantities ${{\hmu}_i}(\la)$:
\be
{\hmu}_i(\la)=\frac{1- q^{-2\, (\la_i+m-i)}}{q-\qq} = q^{-(\lambda_i+m-i)}(\lambda_i+m-i)_q,\quad 1\leq i \leq m.
\label{hmul}
\ee

If  $R\to P$ as $q\to 1$, then by passing to the limit $q\rightarrow 1$ in the above relation we obtain  the quantities ${\hmu}_i(\la)$ for the matrix $\hL$ generating the algebra
$U(gl(N))$:
\be
\hmu_i(\la)=\la_i+m-i,\quad 1\leq i\leq m.
\label{PP}
\ee
Thus, we recover the formula \cite{PP}, expressing the power sums $p_k(\hL)$ via the components $\la_i$ (note that in this case   $m=N$).

Now, we go back to the power sums. As we noted in the previous section, the power sum $p_k(L)$ in the RE algebra $\LL(R)$ corresponds to the element
$z=\tau_{k-1}\, \tau_{k-1}\dots \tau_{1}\in \H(q)$ with respect to the characteristic map and can be written in a simple form $p_k(L) = \Tr_R(L^k)$. Therefore, it is natural to
introduce the power sums in the algebra $\hLL(R)$ by a similar formula $p_k(\hL)=\Trr \hL^k$.

The power sums $p_k(L)$ in the algebra $\LL(R)$ are parameterized in terms of the quantum eigenvalues $\mu_i$ as follows (see \cite{GPS4}):
\be
p_k(L)=\sum_{i=1}^m\, \mu_i^k\,d_i\,,\qquad d_i=q^{-1}\, \prod_{p\not= i}^m \, \frac{\mu_i-q^{-2}\, \mu_p}{ \mu_i-\mu_p}.
\label{pow-sum}
\ee

As for the  power sums $p_k(\hL)$, their parametrization is as follows
\be
p_k(\hL)=\sum_{i=1}^m\, \hmu_i^k\,\hd_i\,,\qquad
\hd_i=q^{-1}\, \prod_{p\not= i}^m \, \frac{\hmu_i-q^{-2}\, \hmu_p-\qq}{ \hmu_i-\hmu_p}.
\label{PPP}
\ee

It should be emphasized that the quantities $p_k(L)$ and $p_k(\hL)$ are {\it symmetric polynomials} in the eigenvalues $\mu_i$ and $\hmu_i$ respectively.

By substituting in (\ref{pow-sum}) and (\ref{PPP}) the characters $\mu_i(\la)$ and $\hmu_i(\la)$ from (\ref{mul}) and (\ref{hmul}), we get
functions in $\la_i$, which are not polynomials.  Also, they are not symmetric in $\la_i$ any more. We consider them to be quantum analogs of shifted symmetric polynomials.
Thus, we have two versions of such polynomials, one for the power sums $p_k(L)$ and another one for $p_k(\hL)$. Observe that our version of quantum shifted symmetric
functions is motivated by the representation theory of the extended algebras $\LL(R)$ and $\hLL(R)$.

\section{Quantum weight systems}

 According to the  approach from \cite{Y1}, we consider the weight systems associated with elements of the symmetric group. More precisely, we consider  the maps
 $$
{w}_{gl(N)}: \SS_n\to U(gl(N)),\qquad \si\mapsto w_{gl(N)}(\si),\quad n\ge 1,
$$
where     $\sigma \in \SS_n$, is a permutation and the quantity $w_{gl(N)}(\si)$ is defined as follows
\be
w_{gl(N)}(\si)= \sum_{i_1\dots i_n=1}^N\,\, \hl_{i_1}^{\, i_{\si(1)}}\, \hl_{i_2}^{\, i_{\si(2)}}\dots \hl_{i_n}^{ \,i_{\si(n)}}.
\label{dva}
\ee
Here $\hl_{i }^{j}$ are the generators of the algebra $gl(N)$ (see the Introduction of the paper). It is known that the element $w_{gl(N)}(\si)$ belongs to the center of the algebra $U(gl(N))$.
Thus, we have a map $\SS\to Z(U(gl(N)))$. Usually, this map is composed with the Harish-Chandra isomorphism, sending the center $Z(U(gl(N)))$ to the algebra $U(\hhh)$.

Our immediate objective is to present the element (\ref{dva}) in another form, more convenient for the subsequent generalization. Let $\si_k$, $1\le  k\le n-1$ be the set of
transpositions\footnote{They are subject to the relations (\ref{H-alg}) at $q=1$.} generating the group $\SS_n$. Since any element $\si\in \SS_n$ can be presented as a
monomial in the generators $\si_k$, then the image of $\sigma$ in the $R$-matrix representation with $R=P$ (i.e.  the usual flip) is a monomial
$\rho_P(\sigma) = \PP_\si(P_1\dots P_{n-1})$. Observe that we endow the set $\{1\dots n\}$ with the left action of the symmetric group.

\begin{proposition}\label{prop:3}
Let $\hL=\|\hl_i^j\|$ be the generating matrix of the algebra $U(gl(N))$. The polynomial {\rm (\ref{dva})} can be written as follows:\rm
\be
w_{gl(N)}(\si)=\Tr_{(1\dots n)} \hL_{\ov{1\to n}}\,\PP_\si(P_1, P_2\dots P_{n-1})
\label{tri}
\ee\it
Here $\hL_{\ov{1\to n}} = \hL_1\hL_2\dots \hL_n$, since for $R=P$ one gets $\hL_{\overline k} = \hL_k$ $\forall\,k\ge 2$.
\end{proposition}

{\bf Proof.} This claim  follows immediately from the fact that the matrix of the usual flip is $P_{i_1 i_2}^{j_1 j_2}=\de_{i_1}^{j_2}\, \de_{i_2}^{j_1}$.
\hfill\rule{6.5pt}{6.5pt}

\medskip

The way of defining the weight systems exhibited in the Proposition \ref{prop:3} enables us to  generalize the notion of the weight system to the algebras $\LL(R)$ and $\hLL(R)$.
Let $z$ be an element of the Hecke algebra $H_n(q)$. It is a polynomial in the generators $\tau_i$: $z=\PP(\tau_1,\dots,\tau_{n-1})$. We define the corresponding quantum
weight system as follows:
\be
{w_{\LL(R)}(z)} = {\Tr_{R(1\dots n)}}{L_{\ov{1\to n}}}\, \PP(R_1,\dots,R_{n-1})\quad \mathrm{or} \quad
{w_{\hLL(R)}(z)}= {\Tr_{R(1\dots n)}}\hL_{\ov{1\to n}}\, \PP(R_1,\dots,R_{n-1}),
\label{zz}
\ee
depending on the algebra which we are dealing with.
According to the result of \cite{IP} the element $w_{\LL(R)}(z)$ belongs to the center of the algebra $\LL(R)$.

\begin{remark} \rm
It should be emphasized that the claim from \cite{IP} and the method of its proof are still valid if a braiding $R$ is a skew-invertible involutive symmetry (i.e. $R^2=I$).
Together with Proposition \ref{prop:3} this claim entails that the elements defined in (\ref{dva}) are central.
\end{remark}

Applying to the element ${w_{\LL(R)}(z)}$ defined in (\ref{zz}) the quantum Harish-Chandra morphism, we get a symmetric polynomial in $\mu_i$. Computing this polynomial is a
somewhat complicated task even in the classical setting. We constrain ourselves to a few examples.

First, consider the element $z=\tau_{n-1}\tau_{n-2}\dots \tau_{2}\tau_{1}\in H_n(q)$. Upon computing the corresponding elements defined in (\ref{zz}), we
get
$$
w_{\LL(R)}(z)=\Tr_{R} L^n,\quad w_{\hLL(R)}(z)=\Tr_{R} \hL^n.
$$
The proof of these formulae is straightforward and is left to the reader. These and all subsequent  formulae can be established by means of the relations
\be
\Tr_{R(k+1)} R_{k}=I,\qquad  \Tr_{R(k+1)} \tilde{L}_{\ov{k+1}}=\Tr_{R(k)} \tilde{L}_{\ov{k}},
\label{zzz}
\ee
where $\tilde{L}$ stands for  $L$ or $\hL$.

\smallskip

\noindent
 Let $n=3$ and take the element $z=\tau_1\tau_2\in {\Bbb H}_3(q)$. Then the corresponding element $w_{\LL(R)}(z)$ is
$$
w_{\LL(R)}(z)={\Tr_{R(123)}}\, {L_{\ov{1\to 3}}}\, R_1R_2=\Tr_{R(123)}\, L_1\, R_1\, L_1\, \RR_1\, R_2\, R_1\, L_1\, \RR_1\, \RR_2\, R_1\, R_2.
$$
The product of the last four $R$-matrices can be simplified with the use of braid relation:
$$
\RR_1\,\underline{ \RR_2\, R_1\, R_2 }= \RR_1\underline{R_1R_2 \RR_1 }= R_2 \RR_1
$$
and we come to
$$
w_{\LL(R)}(z) = \Tr_{R(123)}\, L_1 R_1 L_1 \RR_1 R_2 R_1R_2  L_1\RR_1.
$$

Then we again use the braid relation for the product of four $R$-matrices and compute the $R$-trace in the 3-rd space according to the first formula of (\ref{zzz}):
$$
w_{\LL(R)}(z) = \Tr_{R(12)}\,L_1 \underline{R_1 L_1 R_1 L_1}\RR_1 = \Tr_{R(12)}\, L_1^2 R_1 L_1 = \Tr_R L^3.
$$
Here for the underlined product we apply the defining relations (\ref{REA}) of RE algebra $\LL(R)$. Thus, this quantum weight system does not distinguish the elements
$z=\tau_1\tau_2$ and $z=\tau_2\tau_1$. In general, it is possible to show that the weight system associated with the algebra $\LL(R)$ does not distinguish the elements
$z=\tau_1\tau_2\dots \tau_k$ and $z=\tau_k\dots \tau_2\tau_1$. But it is not so for the weight system  $w_{\hLL(R)}(z)$ associated with the algebras $\hLL(R)$.

Indeed, with the same calculations as above we arrive at the expression:
$$
w_{\hLL(R)}(z)={\Tr_{R(123)}}\, {\hL_{\ov{1\to 3}}}\, R_1R_2 = \Tr_{R(12)}\, \hL_1\, \underline{R_1 \hL_1  R_1 \hL_1}  \RR_1.
$$
But now for the transformation of the underlined product we have to use the defining relations (\ref{odin}) of the algebra $\hLL(R)$ and obtain the following result:
$$
 w_{\hLL(R)}(z)=\Tr_{R(12)}\, (\hL_1^2 R_1\hL_1 + \hL_1 R_1 \hL_1 \RR_1-\hL_1^2) = \Trr\hL^3+(\Trr \hL)^2- \Trr\hL^2.
$$
This result looks like the classical one. The difference becomes transparent after applying the quantum Harish-Chandra morphism
$$
w_{\hLL(R)}(z)=\sum_i (\hmu_i^3+\hmu_i^2)\, \hd_i+(\sum_i \hmu_i\, \hd_i)^2,
$$
where $\hd_i$ are defined in (\ref{PPP}).

Consider one more example for $n=3$. We choose  $z=\tau_1\tau_2\tau_1=\tau_2\tau_1\tau_2$. Then we have:
\begin{eqnarray*}
w_{\LL(R)}(z) &=&\Tr_{R(123)} (L_{\overline 1}L_{\overline 2}L_{\overline 3}R_2 R_1 R_2)=\Tr_{R(123)} (L_{1}L_{\overline 2}R_2 R_1 L_{1}R_2) \\
 &=& \Tr_{R(123)}( L_{1}R_1L_{1}R_2 R_1 L_{1}) =\Tr_{R(12)} (L_{1}R_1L_{1} R_1L_{1})  \\
 &=& \Tr_{R(12)} (L_{1}R_1L_{1} (R_1^{-1} +(q-q^{-1})I)L_{1}) = \Trr( L)\,\Trr(L^2)+(q-\qq)\Trr(L^3). \\
\end{eqnarray*}

Note that we did not use here the defining relations of the algebra, so we get the same result of computations for the corresponding element
$w_{\hLL(R)}(z)\in Z(\hLL(R))$.

\section{General type algebras and concluding remarks}

So far we have been dealing with even symmetries of bi-rank  $(m|0)$. As for the algebras related to the symmetries of bi-rank $(m|n)$, the Cayley-Hamilton polynomial for their generating
matrices $L$ was obtained in \cite{GPS1}. Its coefficients  are also central and (upon multiplying by a central element) the polynomial can be factorized
in a product of two polynomials with central coefficients. This property enables us to introduce two families of quantities $\{\mu_i\}$, $1\leq i\leq m$ and  $\{\nu_i\}$, $1\leq i\leq n$,
which are the roots of these two polynomial factors respectively. In this case the parameterization of the power sums $p_k(L)$ in terms of quantum eigenvalues takes  the following form
\cite{GPS2}:
\be
p_k(L)=\sum_{i=1}^m\, \mu_i^k\,d_i+\sum_{j=1}^n\, \nu_j^k \, {f_j}, \label{pk} \ee
where
\be
d_i=q^{-1}\, \prod_{p\not= i}^m \, \frac{\mu_i-q^{-2}\, \mu_p}{ \mu_i-\mu_p}\, \prod_{j=1}^n \frac{\mu_i-q^{2}\, \nu_j}{ \mu_i-\nu_j},\qquad
{f_j}=-q\, \prod_{p\not=j}^n \frac{\nu_j-q^{2}\, \nu_p}{ \nu_j-\nu_p}\, \prod_{i= 1}^m \, \frac{\nu_j-q^{-2}\, \mu_i}{ \nu_j-\mu_i}. \label{df} \ee

Observe that the power sums $p_k(L)$ in this case are the super-symmetric polynomials in the quantities $\qq \mu_i$ and $q \nu_i$ as defined in \cite{S}.

Again, we treat  the map, sending the  elements $p_k(L)$ into the above super-symmetric polynomials in eigenvalues $\mu_i$ and $\nu_i$, as a quantum analog of
the Harish-Chandra morphism in the RE algebras of general type. In order to get a similar morphism in the modified RE algebras it suffices to replace
$\mu_i$, $\nu_i$, $d_i$ and $f_i$ in the formulae (\ref{pk}) for respectively $\hmu_i$, $\hat{\nu}_i$, $\hd_i$ and $\hf_i$, where the  quantities $\hd_i$ and $\hf_i$
can be deduced from (\ref{df}) by means of the relations (\ref{zamena}) and similar formulae for $\nu_i$.

However,  the characters of the spectral values $\nu_i$ in the modules $\vl$ are not yet computed: so far we do not know the relations between the quantities $\nu_i$ and $\la_i$.

Summing up, we want to emphasize the universal character of our method: it is valid in the RE algebras (modified or not), associated with any skew-invertible
Hecke symmetry. We do not use any object of a quantum group type, instead we  use the eigenvalues of the generating matrices arising from the corresponding
Cayley-Hamilton identities.

In this connection we want to mention the paper \cite{JLM}, where the characters of the power sums (called in that paper the quantum Gelfand invariants) are computed
in the particular case related to the quantum group $U_q(gl(N))$. The main tool in this way is the RTT presentations of $U_q(gl(N))$ and its affine analog going back to the
Leningrad school \cite{KS, RTF}. In fact, this construction enables one to map the RE algebra, corresponding to the Hecke symmetry $R$ coming from  $U_q(gl(N))$, into this
quantum group. This map enables one to realize entries of the  generating matrix of the RE algebra as elements from $U_q(gl(N))$. Thus, one gets the mentioned in footnote 2
matrix which meets  a Cayley-Hamilton identity with central coefficients.

Note that the characters of the power sums $\chi_\la(p_k(L))$ computed in \cite{GSZ} coincide with those obtained in \cite{JLM} and they do not depend on
the Hecke symmetry $R$. The only difference is that the quantity  $N=\dim V$ entering the mentioned characters should be replaced by the degree of the Cayley-Hamilton polynomial.
It is not a surprising property since usually  numerical quantities arising in $R$-matrix technique can be expressed via the bi-rank of the  initial Hecke symmetry.

\end{document}